\input amstex 
\documentstyle{amsppt} 
\input xy
\xyoption{all}
\loadbold
\magnification=1200
\pagewidth{6.4truein}
\pageheight{8.87truein}

\redefine\cdot{\boldsymbol\cdot}
\define\Hom{\operatorname{Hom}}
\define\colim{\operatornamewithlimits{colim}}

\NoRunningHeads

\topmatter 

\title   Model Structures and the Oka Principle  \endtitle
\author  Finnur L\acuteaccent arusson  \endauthor
\affil   University of Western Ontario \endaffil
\address Department of Mathematics, University of Western Ontario,
         London, Ontario N6A~5B7, Canada \endaddress
\email   larusson\@uwo.ca \endemail

\thanks The author was supported in part by the Natural Sciences and 
Engineering Research Council of Canada. \endthanks

\thanks First version March 2003; revised September 2003; minor
changes 2 February 2004. \endthanks

\subjclass Primary: 32Q28; secondary: 18F10, 18F20, 18G30, 18G55, 
32E10, 55U35 \endsubjclass

\abstract  We embed the category of complex manifolds into the 
simplicial category of prestacks on the simplicial site of Stein 
manifolds, a prestack being a contravariant simplicial functor from the
site to the category of simplicial sets.  The category of prestacks 
carries model structures, one of them defined for the first time here, 
which allow us to develop {\it holomorphic homotopy theory}.  More 
specifically, we use homotopical algebra to study lifting and extension
properties of holomorphic maps, such as those given by the Oka Principle.
We prove that holomorphic maps satisfy certain versions of the Oka
Principle if and only if they are fibrations in suitable model
structures.  We are naturally led to a simplicial, rather than a
topological, approach, which is a novelty in analysis.
\endabstract

\endtopmatter

\document

\subhead 1. Introduction  \endsubhead  
This paper, like its predecessor \cite{L}, is about model structures
in complex analysis.  Model structures are good for many things, but
here we view them primarily as a tool for studying lifting and extension
properties of holomorphic maps, such as those given by the Oka Principle.
More precisely, model structures provide a framework for investigating 
two classes of holomorphic maps such that the first has the right lifting 
property with respect to the second and the second has the left lifting
property with respect to the first in the absence of topological
obstructions.   (It is more natural, actually, to consider homotopy
lifting properties rather than plain lifting properties.)  We seek to 
make the maps in the first class into fibrations and those in the second
class into cofibrations, with weak equivalences being understood in the
topological sense.  The machinery of abstract homotopy theory can then be
applied.

The version of the Oka Principle we focus on here involves the inclusion
$T\to S$ into a Stein manifold of a closed complex submanifold and a
holomorphic fibre bundle $X\to Y$ whose fibre is an elliptic manifold.
Loosely speaking, ellipticity means receiving many holomorphic maps from
Euclidean spaces; it is thus dual to being Stein.  A deep theorem of
Gromov \cite{G, FP} implies that for any commuting square
$$\xymatrix{
T \ar[r] \ar[d] & X \ar[d] \\
S \ar[r] & Y }$$
in which $T\to X$ and $S\to Y$ are otherwise arbitrary holomorphic
maps, the inclusion of the space of holomorphic liftings $S\to X$ into
the space of continuous liftings is a weak equivalence in
the compact-open topology.  Since $T\to S$ is a topological cofibration
and $X\to Y$ is a topological fibration, it follows by basic topology
that there is a holomorphic lifting if one of the two maps is a homotopy
equivalence.  This looks very much like a holomorphic manifestation of
Quillen's first axiom for a model category \cite{Q, Ch\. I, p\. 0.1}, so
it is natural to ask whether there is a model category containing the
category of complex manifolds in which Stein inclusions are cofibrations,
weak equivalences are defined topologically, and being a fibration is
equivalent to an Oka property, such as the one attributed to elliptic
bundles by Gromov's theorem.  The main result of this paper is that a
stronger and perhaps more natural Oka property, in which we consider
not a single square but a continuous family of them, is equivalent to
fibrancy in a new model category containing the category of complex
manifolds.  Elliptic manifolds are fibrant in this new sense, but it is
still an open question whether all elliptic bundles are fibrations.

We equip the category of Stein manifolds in a natural way with a 
simplicial structure and a compatible topology, turning it into a
simplicial site, and embed the category of complex manifolds into the 
simplicial category of prestacks on this site.  By a prestack we 
mean a contravariant simplicial functor from the site to the category of
simplicial sets.  We make use of recent work of To\"en and Vezzosi 
\cite{TV}, generalizing the homotopy theory of simplicial presheaves on 
ordinary, discrete sites to prestacks on simplicial sites.  The category of
prestacks carries several interesting model structures.  Strengthening the 
main result of \cite{L}, we show that the prestack represented by a complex 
manifold $X$ is fibrant in the so-called projective structure (so $X$
represents a stack, in the terminology of \cite{TV}) if and only if $X$ 
satisfies what we call the weak Oka property.  This means that for every Stein
manifold $S$, the inclusion of the space of holomorphic maps from $S$ to $X$
into the space of continuous maps is a weak equivalence in the compact-open
topology.  By Gromov's theorem, this holds if $X$ is elliptic.
We generalize the weak Oka property to holomorphic maps (viewing manifolds
as constant maps) and show that it is equivalent to being a projective
fibration.  

We introduce a new simplicial model structure on the category of prestacks on
the Stein site, in a sense the smallest one in which every Stein inclusion is
a cofibration.  We characterize the fibrations in this structure and show
that a holomorphic map is a fibration if and only if it satisfies a new,
stronger Oka property.  This Oka property is defined explicitly in purely
analytic terms, without reference to, but with guidance from, abstract 
homotopy theory.  For a holomorphic map which is a homotopy equivalence,
it turns out to be simply the homotopy right lifting property with respect
to all Stein inclusions.  By Gromov's theorem, elliptic manifolds are fibrant.
I conjecture that this extends to nonconstant maps: that elliptic bundles
are fibrations.  So far, this is known for covering maps but remains open
for nontrivial bundles in general.

The interface between  complex analysis and homotopical algebra will be
explored further in future work.  For more motivation, see the final
remarks at the end of the paper, and for more background, the
introduction in \cite{L} and the survey \cite{F2}.

\smallskip\noindent
{\it Acknowledgement.}  I am indebted to Rick Jardine for helpful
conversations.

\subhead 2. The embedding \endsubhead
Let $\Cal M$ be the category of complex manifolds, second countable but
not necessarily connected, and holomorphic maps.  As the first step in 
the development of holomorphic homotopy theory, or more specifically a
homotopy-theoretic study of the Oka Principle, we wish to embed $\Cal M$
in a simplicial model category.

Now $\Cal M$ has a natural simplicial structure (enrichment over the 
category $s\bold{Set}$ of simplicial sets), making it a simplicial
object in the category of categories with a discrete simplicial class of
objects.  For complex manifolds $X$ and $Y$, the mapping space 
$\Hom(X,Y)$ is the singular set $s\Cal O(X,Y)$ of the space of 
holomorphic maps from $X$ to $Y$ with the compact-open topology.

Let $\Cal S$ be the full subcategory of Stein manifolds with this 
simplicial structure.  It is a small category, or at least equivalent
to one, since a connected Stein manifold can be embedded into Euclidean
space.  A prestack on $\Cal S$ (in the terminology of \cite{TV}) is a
contravariant simplicial functor (morphism of simplicial categories) 
$\Cal S\to s\bold{Set}$.  Let $\frak S$ denote the category of 
prestacks on $\Cal S$ with its own natural simplicial structure (in a
sense that is stronger than the sense in which $\Cal M$ is a simplicial
category; see \cite{GJ, IX.1}).

By the simplicial Yoneda lemma \cite{GJ, IX.1.2}, if $S$ is an
object of $\Cal S$ and $F$ is a prestack on $\Cal S$, then there is a
natural isomorphism of simplicial sets
$$F(S)\cong \Hom_{\frak S}(\Hom_{\Cal M}(\cdot,S),F).$$
(From now on we will usually omit the subscripts.)  Hence there is a 
simplicially full embedding of $\Cal S$ into $\frak S$, taking an 
object $S$ of $\Cal S$ to the prestack $\Hom(\cdot,S)$ represented by 
$S$.

The embedding $\Cal S\to\frak S$ clearly extends to a functor $\Cal
M\to\frak S$, taking a complex manifold $X$ to the prestack
$\Hom(\cdot,X)$ on $\Cal S$ represented by $X$.  This functor induces
monomorphisms (injections at each level) of mapping spaces, as is
easily seen by plugging in the terminal object of $\Cal S$, the 
one-point manifold $\frak p$.  Hence, for complex manifolds $X$ and
$Y$, we have a monomorphism of mapping spaces
$$\Hom(X,Y) \to \Hom(\Hom(\cdot,X),\Hom(\cdot,Y)),$$
which is an isomorphism when $X$ is Stein, and we have a simplicial
embedding of $\Cal M$ into $\frak S$.  Whether the embedding is full
remains to be investigated.

\subhead 3. Remarks \endsubhead  
We would like to motivate the above construction and explain why it seems
to produce an appropriate setting for applying homotopical algebra
in complex analysis.  Yoneda embeddings provide the canonical way 
of closing geometric categories under limits and colimits.  This is the 
first step in the homotopy theory of schemes, for instance; I know of no 
alternative.  In our paper \cite{L}, we embedded $\Cal M$ into the category
of all simplicial presheaves on $\Cal S$, but there is every reason to take 
into account the topology on our hom-sets and restrict attention to those
simplicial presheaves that respect it, now that the homotopy theory of
simplicial presheaves on ordinary, discrete sites has been generalized to
prestacks on simplicial sites by To\"en and Vezzosi \cite{TV}.  Indeed, 
we want a full embedding of the category of complex manifolds into a 
simplicial model category, at least for Stein sources, and with plain
simplicial presheaves we cannot expect this.  Homotopy theory gives
information about simplicial hom-sets and maps between them; to apply 
such results in complex analysis, we need to know that simplicial 
hom-sets essentially equal spaces of holomorphic maps.  We get this at 
least when the source is Stein; this has proved sufficient so far.

It would seem simpler and more natural to use presheaves of topological 
spaces on $\Cal S$ rather than simplicial presheaves.  The homotopy 
theory of the former is not available in the literature --- although it
could presumably be developed in a straightforward manner for a suitable
locally presentable category of topological spaces, now that one such has
been discovered: J\. Smith's category of $I$-spaces --- but that is not
why we use the latter.  The reason is that we are aiming for a model
structure in which the inclusion $T\hookrightarrow S$ of a closed complex
submanifold $T$ in a Stein manifold $S$ is a cofibration (this is the
intermediate structure,
defined below).  It is appropriate, then, to require such an inclusion to
induce a pointwise cofibration, so in the topological setting we would 
need $\Cal O(X,T)\to\Cal O(X,S)$ to be a cofibration of topological spaces 
for every Stein manifold $X$.  There are simple examples for which this
fails.  For instance, let $S$ be the complex plane with a puncture, $T$ be
a one-point subset of $S$, and $X$ be the complex plane with the integers 
removed.  Then $\Cal O(X,T)\to\Cal O(X,S)$ is not a cofibration, not
even in the weaker of the two senses considered by topologists, because
the point $\Cal O(X,T)$ in the space $A=\Cal O(X,S)$ does not have a 
neighbourhood contractible in $A$.  Indeed, there are uncountably many 
homotopy classes of holomorphic maps $X\to S$ (consider winding numbers 
around each integer), so $A$ has uncountably many connected components, 
and every nonempty open subset of $A$ contains uncountably many of these,
so it is not contractible in $A$.  However, the induced map $s\Cal O(X,T)
\to s\Cal O(X,S)$ is a cofibration of simplicial sets, simply because it
is injective at each level.  Shifting our focus from the spaces of
holomorphic maps themselves to the singular sets that catalogue continuous
families of holomorphic maps with nice parameter spaces alleviates the
difficulties associated with the compact-open topology for noncompact
sources.

Thus we are, somewhat surprisingly, led to a simplicial approach, which
is a novelty in analysis.  Fortunately, there is often no loss involved 
in applying the singular functor to spaces of holomorphic maps, because
the singular functor not only preserves but also reflects fibrations.  
For example, if $A$ and $B$ are spaces of holomorphic maps
and $A\to B$ is a map such that the induced map $sA\to sB$ of
mapping spaces is a Kan fibration, as might follow from some 
homotopy-theoretic arguments, then $A\to B$ itself is a Serre fibration
(and conversely).  Also, $sA\to sB$ is a weak equivalence if and only
if $A\to B$ is.

\subhead 4. The projective model structures  \endsubhead  
The category $\frak S$ carries several interesting simplicial model 
structures.  We begin by describing the most basic one, the coarse
projective structure, originally defined by Dwyer and Kan \cite{GJ, 
IX.1}.  (We call it coarse because it is associated to the coarsest 
topology on $\Cal S$, that is, the trivial topology; see below.)
In this structure, which is cofibrantly generated and proper, weak 
equivalences and fibrations are defined pointwise, so a map $F\to G$ of
prestacks on $\Cal S$ is a weak equivalence or a fibration if the 
component maps $F(S)\to G(S)$ are weak equivalences or fibrations of 
simplicial sets, respectively, for all objects $S$ in $\Cal S$.  In 
particular, a holomorphic map $X\to Y$, viewed as a map of the 
prestacks represented by $X$ and $Y$, is a weak equivalence or a 
fibration in the coarse projective structure if the induced maps 
$\Cal O(S,X)\to\Cal O(S,Y)$ are weak equivalences or Serre fibrations
of topological spaces, respectively, for all Stein manifolds $S$.  
Cofibrations are defined by a left lifting property.  The prestacks
represented by Stein manifolds are both cofibrant and fibrant.

Now we move to the projective structure on $\frak S$, which is obtained
by a left Bousfield localization of the coarse projective structure.
There will be a larger class of weak equivalences, defined using a topology 
on the simplicial category $\Cal S$, turning it into a simplicial site.
The cofibrations are the same as in the coarse projective structure, so 
they can be referred to simply as projective cofibrations.  The projective 
fibrations are determined by a right lifting property; they form a
subclass of the class of pointwise fibrations.

The category of components $c\Cal S$ (also called, at some 
risk of confusion, the homotopy category) of the simplicial category 
$\Cal S$ has the same objects as $\Cal S$, and its hom-sets are the sets
of path components of the simplicial hom-sets of $\Cal S$.  We can also 
obtain $c\Cal S$ from $\Cal S$ by identifying maps in the underlying 
category of $\Cal S$ that can be joined by a string of homotopies
(provided by the simplicial structure).  By precomposition by the morphism 
$\Cal S\to c\Cal S$, a presheaf on $c\Cal S$ gives a presheaf on $\Cal S$ 
such that equivalent maps in $\Cal S$ induce the same restriction maps.  
Conversely, such a presheaf on $\Cal S$ descends to $c\Cal S$.  Prestacks 
respect the simplicial structure, so they preserve homotopies, so the 
homotopy presheaves of a prestack on $\Cal S$ naturally live on $c\Cal S$
(or, more precisely, on overcategories thereof).

A topology on $\Cal S$, turning it into a simplicial site (an $S$-site in 
the language of \cite{TV}), is a Grothendieck topology in the usual sense on
the category of components $c\Cal S$.  A map of prestacks is a weak 
equivalence, or acyclic, with respect to the topology, if it induces 
isomorphisms of homotopy sheaves in all degrees, that is, isomorphisms of 
the sheafifications (with respect to the given topology) of homotopy 
presheaves in all degrees.  By a theorem of To\"en and Vezzosi \cite{TV, 
Thm\. 3.4.1}, the projective structure on $\frak S$ is a cofibrantly 
generated, proper, simplicial model structure. 

The projective structure specializes in two ways.  It equals the coarse 
projective structure when the topology on $c\Cal S$ is trivial.  Also, 
when the simplicial structure on $\Cal S$ is trivial (discrete), so $c\Cal
S=\Cal S$, then $\Cal S$ is an ordinary site and we obtain the well-known
projective structure (sometimes called local) for simplicial presheaves 
on $\Cal S$.

The topology we shall put on the Stein site $\Cal S$ is the \lq\lq 
usual\rq\rq\ topology employed in \cite{L}, except we now view it as a 
topology on the category of components $c\Cal S$, which is obtained from 
the plain category of Stein manifolds and holomorphic maps by identifying
holomorphic maps $X\to Y$ that are homotopic in the usual sense that 
they can be joined by a continuous path in $\Cal O(X,Y)$ with the
compact-open topology.  In other words,
$$\hom_{c\Cal S}(X,Y)=\pi_0\Cal O(X,Y).$$
A cover of a Stein manifold $S$ is a family of holomorphic maps into 
$S$ such that by suitably deforming each map $X\to S$ inside $\Cal
O(X,S)$, we get a family of biholomorphisms onto Stein open subsets of
$S$ which cover $S$.  This defines a Grothendieck topology on $c\Cal S$.

The acyclic maps have a very simple description.  First, for any map
from the point $\frak p$ to an open ball $B$, the map $\frak p\to B\to 
\frak p$ is the identity and the map $B\to\frak p\to B$ is homotopic to 
the identity through holomorphic maps keeping the image point of the map
$\frak p\to B$ fixed.  Hence, if $F$ is a prestack on $\Cal S$, the 
restriction map $F(\frak p)\to F(B)$ is a homotopy equivalence, in fact
the inclusion of a strong deformation retract.  Since every cover has a
refinement by balls, this implies that a map $F\to G$ of prestacks on 
$\Cal S$ is acyclic if and only if $F(\frak p)\to G(\frak p)$ is acyclic.
Here it is crucial that prestacks respect the simplicial structure on 
$\Cal S$; this does not work for arbitrary simplicial presheaves.  It 
follows that a holomorphic map $f:X\to Y$ of complex manifolds, viewed as
a map of the prestacks represented by $X$ and $Y$, is acyclic if and only
if it is a topological weak equivalence, that is, a homotopy equivalence. 

\subhead 5. The injective model structures  \endsubhead
We will also need the so-called injective model structures on $\frak S$
\cite{TV, 3.6}.  The coarse injective structure is a proper, simplicial
model structure on $\frak S$ in which weak equivalences and cofibrations
are defined pointwise and fibrations are defined by a right lifting
property.  In the injective structure, which is also proper and 
simplicial, the cofibrations are the same, weak equivalences are acyclic
with respect to the chosen topology on $\Cal S$, and fibrations are
defined by a right lifting property.  Injective cofibrations are and 
will be referred to simply as monomorphisms.

\subhead 6. A Quillen equivalence \endsubhead
Consider the functor $P:\frak S\to\frak S$ taking a prestack $F$ to the
prestack $PF=\Hom(s\cdot, F(\frak p))$ and taking a map $f:F\to G$ to the 
map $Pf:PF\to PG$ induced by the map $F(\frak p)\to G(\frak p)$.  This 
functor is a projection: $P\circ P=P$.  There is a natural transformation
$\eta$ from the identity functor on $\frak S$ to $P$:  if $F$ is a prestack 
and $S$ is an object of $\Cal S$, the map (morphism of simplicial sets)
$\eta_F(S): F(S)\to PF(S)= \Hom(sS,F(\frak p))$ comes from the map 
$sS=\Hom(\frak p,S)\to\Hom(F(S), F(\frak p))$ given directly by $F$.  
Here, again, it is crucial that prestacks respect the simplicial structure
on $\Cal S$; this does not work for arbitrary simplicial presheaves.  The
square
$$\xymatrix{
F \ar[r]^f \ar[d]_{\eta_F} & G \ar[d]^{\eta_G} \\
PF \ar[r]^{Pf} & PG
}$$
commutes simply because maps of prestacks commute with restrictions.  Note
that the map $\eta_{PF}=P(\eta_F): PF\to P^2F=PF$ is the identity.  Also, 
$\eta_F:F \to PF$ is acyclic, since $\eta_F(\frak p)$ is the identity.  
The pair $P$, $\eta$ is a key element of the structure of $\frak S$ and 
plays an important role in our theory.  It is an example of what is
called a localization functor.

If $A$ is a simplicial set, let $\tilde A$ denote the constant prestack
with $\tilde A(S)=A$ for each $S$ in $\Cal S$ and with all restriction 
maps equal to the identity.  Define a functor $R:s\bold{Set}\to\frak S$
by $RA=P\tilde A=\Hom(s\cdot, A)$.  A map $f$ from a prestack $F$ to $RA$
factors as
$$\xymatrix{
F \ar[r]^f \ar[d]_{\eta_F} & RA \ar@{=}[d]^{\eta_{RA}} \\
PF \ar[r]^{Pf} & RA 
}$$
so $f$ is determined by $Pf$, which is induced by the map $F(\frak p)\to
RA(\frak p)=A$.  Hence, we have a pair of adjoint functors
$$L:\frak S\to s\bold{Set}:R, \qquad LF=F(\frak p), \qquad
RA=\Hom(s\cdot, A),$$
with a natural bijection
$$\hom_{\frak S}(F,RA) \cong \hom_{s\bold{Set}}(LF,A)$$
for every prestack $F$ on $\Cal S$ and every simplicial set $A$.

We see that a map $F\to RA$ is acyclic if and only if the corresponding
map $LF\to A$ is.  Also, it is clear that $L$ takes monomorphisms to 
cofibrations and preserves weak equivalences.
Hence, $(L,R)$ is a pair of Quillen equivalences between 
the category of simplicial sets and the category of prestacks on $\Cal
S$ with the projective structure or the injective structure \cite{H, 8.5}.
Such a pair induces equivalences of homotopy categories, so the homotopy
category of $\frak S$ is the ordinary homotopy category of simplical
sets or topological spaces.  It also follows that $R$ takes fibrations
of simplicial sets to injective fibrations; in particular, if $K$ is a 
fibrant simplicial set (a Kan complex), then the prestack $\Hom(s\cdot,
K)$ is injectively fibrant.  Hence, if $X$ is a complex manifold, so 
$\eta_X:X\to PX$ is a monomorphism, then $\eta_X$ is an injectively 
cofibrant fibrant model for $X$.

\subhead 7. Projective fibrations  \endsubhead
The projective structure is the left Bousfield localization of the
coarse projective structure on $\frak S$ with respect to the class of
acyclic maps of prestacks.  The theory of the left Bousfield
localization provides a useful characterization of projective
fibrations.

Let $f:F\to G$ be a pointwise fibration of prestacks such that $F(\frak p)$
and $G(\frak p)$ are fibrant, so $PF$ and $PG$ are injectively and hence 
projectively fibrant.  Then the square
$$\xymatrix{
F \ar[r]^f \ar[d]_{\eta_F} & G \ar[d]^{\eta_G} \\
PF \ar[r]^{Pf} & PG
}$$
is a localization of $f$ \cite{H, 3.2.16}.  The map $f$ is a projective
fibration if and only if this square is a homotopy pullback in the
coarse projective structure \cite{H, 3.4.8}.  This means that the
natural map from $F$ to the homotopy pullback of $G\to PG\gets PF$ is
pointwise acyclic.  Since $F(\frak p)\to G(\frak p)$ is a fibration,
$Pf$ is a pointwise fibration, so the homotopy pullback is naturally
pointwise weakly equivalent to the ordinary pullback (taken pointwise).

In summary, a map $F\to G$ of prestacks fibrant at $\frak p$ is a projective
fibration if and only if it is a pointwise fibration and the induced map $F
\to G\times_{PG} PF$ is pointwise acyclic.  In particular, a prestack 
$F$ is projectively fibrant if and only if it is pointwise fibrant and 
$\eta_F$ is pointwise acyclic.

\subhead 8. Stacks on the Stein site and the weak Oka property \endsubhead  
A pointwise fibrant prestack on the simplicial site $\Cal S$ is called,
in the language of \cite{TV}, a stack on $\Cal S$ (with respect to the
chosen topology) if it is projectively fibrant.  Loosely speaking, this
is a \lq\lq homotopy sheaf condition\rq\rq, with the limits in the usual
sheaf condition replaced by homotopy limits.  The sheaf condition is not
really relevant here; indeed, the prestacks and the topology live on 
different categories ($\Cal S$ and $c\Cal S$, respectively), so we will
not be talking about a prestack being a sheaf in the usual sense.

We say that a complex manifold $X$ satisfies the weak Oka property, or 
that $X$ is weakly Oka, if the inclusion map $\Cal O(S,X)\hookrightarrow 
\Cal C(S,X)$ is a weak equivalence for all Stein manifolds $S$, where the
spaces of holomorphic and continuous maps from $S$ to $X$ carry the 
compact-open topology.  The main result of \cite{L} characterizes the
weak Oka property (there called the Oka-Grauert property) in terms of 
excision; the following theorem, using a better model structure, is more
to the point.

\proclaim{9. Theorem}  A complex manifold is weakly Oka if and only if
it represents a stack on the Stein site.
\endproclaim

\demo{Proof}  A prestack $F$ is projectively fibrant if and only if it
is pointwise fibrant and the map $\eta_F:F\to PF$ is pointwise acyclic.
If $F$ is represented by a complex manifold $X$, so it is pointwise 
fibrant, this means that the map from $F(S)=s\Cal O(S,X)$ to $PF(S)=
\Hom(sS,sX)=s\Cal C(|sS|,X)$ is acyclic for every Stein manifold $S$.
Since $PF(S)$ is homotopy equivalent to $s\Cal C(S,X)$, this is nothing
but the weak Oka property.
\qed\enddemo

It is an interesting open question whether the inclusions $\Cal O(S,X)
\hookrightarrow \Cal C(S,X)$ have functorial homotopy inverses when $X$
is weakly Oka.  Since the spaces in question are not known to be cofibrant,
even the existence of pointwise homotopy inverses is not clear \cite{L,
Thm\. 2.2}, but it is in the simplicial setting, so we ask whether the 
pointwise homotopy equivalence $\eta_X:X \to PX$ is in fact a simplicial
homotopy equivalence of prestacks.  This would follow if $\eta_X$ was not
only a monomorphism but actually a projective cofibration \cite{H, 9.6.5},
that is, if $PX$ was a cofibrant fibrant model for $X$ not only in the 
injective structure but also in the projective structure.

\subhead 10. The weak Oka property for maps  \endsubhead
Let us generalize the above discussion from objects to arrows.  We 
say that a holomorphic map $f:X\to Y$ satisfies the weak Oka property,
or that $f$ is weakly Oka, if
\roster
\item the induced map $\Cal O(S,X)\to\Cal O(S,Y)$ is a Serre fibration
and
\item the inclusion $\Cal O(S,X) \hookrightarrow \{h\in\Cal C(S,X) :
f\circ h \in \Cal O(S,Y)\}$ is acyclic
\endroster
for every Stein manifold $S$.  

In particular, if $f:X\to Y$ is weakly Oka, then every continuous map $h$
from a Stein manifold to $X$ such that $f\circ h$ is holomorphic
can be continuously deformed through such maps to a holomorphic map.  
Clearly, a complex manifold $X$ is weakly Oka if and only if the constant
map $X\to\frak p$ is weakly Oka.

\proclaim{11. Theorem}  A holomorphic map is weakly Oka if and only if
it is a projective fibration.
\endproclaim

\demo{Proof}  A holomorphic map $f:X\to Y$ is a projective fibration if
and only if it is a pointwise fibration, meaning that the induced map 
$\Cal O(S,X)\to\Cal O(S,Y)$ is a Serre fibration for every Stein manifold
$S$, and the induced map $X\to Y\times_{PY}PX$ is pointwise acyclic,
which is equivalent to the map
$$\Cal O(S,X)\to \Cal O(S,Y)\times_{\Cal C(S,Y)} \Cal C(S,X)$$
being acyclic for every Stein manifold $S$.  Finally, the space on the 
right is the space of continuous maps $h:S\to X$ such that $f\circ h$ is 
holomorphic.
\qed\enddemo

\subhead 12. The intermediate model structure  \endsubhead
We now introduce a new simplicial model structure on $\frak S$,
in between the projective and injective structures in the 
sense that it has fewer fibrations than the projective structure and
more fibrations than the injective structure; for cofibrations it is
the other way around.  The weak equivalences are the same: the maps that
are acyclic with respect to the chosen topology on $\Cal S$.

By a Stein inclusion we mean the inclusion $T\hookrightarrow S$ of a
closed complex submanifold $T$ in a Stein manifold $S$ (then $T$ is also
Stein).  Let the set $C$ consist of all the monomorphisms
$$S\times\partial\Delta^n\cup_{T\times\partial\Delta^n}
T\times\Delta^n \to S\times\Delta^n,$$
in $\frak S$, where $T\hookrightarrow S$ is a Stein inclusion and $n\geq
0$.  Among these maps are the Stein inclusions $T\hookrightarrow S$
themselves (with $n=0$), as well as the standard generating cofibrations
$S\times\partial \Delta^n \to S\times \Delta^n$ for the projective
structure (with $T=\varnothing$).

Let $\Cal C$ be the saturation of $C$, that is, the smallest class of 
maps in $\frak S$ which contains $C$ and is closed under pushouts, 
retracts, and transfinite compositions.  The maps in $\Cal C$ are called
intermediate cofibrations; they are retracts of transfinite compositions
of pushouts of maps in $C$.  An intermediate fibration is defined to be
a map with the right lifting property with respect to all acyclic
intermediate cofibrations.

The idea of an intermediate structure in which Stein inclusions would be
cofibrations came up in a discussion with Rick Jardine, who subsequently 
showed me how to obtain such a structure and later wrote up a proof in 
\cite{J}, which we follow below.  The argument for a simplicial site is 
the same as for the special case of a discrete site, treated in \cite{J}.
Later, I learned that one can show that the intermediate structure exists
and, moreover, is cofibrantly generated, using a very general argument 
due to T\. Beke and J\. Smith \cite{B, Thm\. 1.7}, based solely on 
$\frak S$ being locally presentable and the class of weak equivalences 
being accessible.  (Cofibrant generation is also contained in a second
version of \cite{J}.)  Unfortunately, the generating set of acyclic
cofibrations produced by this method is too large to be of much
practical use.  

\proclaim{13. Theorem}  There is a proper, simplicial model structure on
$\frak S$, called the intermediate structure, with cofibrations, 
fibrations, and weak equivalences defined as above.
\endproclaim

\demo{Proof} Consider factorization first.  Since $\frak S$ is locally
presentable, a standard small object argument shows that a map $X\to Y$
of prestacks can be factored as $X @>j>> Z @>p>> Y$, where $j$ is in 
$\Cal C$ and $p$ has the right lifting property with respect to every
map in $\Cal C$, so $p$ is an acyclic intermediate fibration (note that
we do not know the converse of this yet).

For the other factorization, we make use of the injective structure to
factor $X\to Y$ as $X @>i>> W @>q>> Y$, where $i$ is an acyclic 
injective cofibration and $q$ is an injective fibration and hence an
intermediate fibration.  Then factor $i$ as above as $X @>j>> Z
@>p>> W$, where $j$ is an intermediate cofibration and $p$ is an
acyclic intermediate fibration.  Then $j$ is acyclic too and $qp$ is an
intermediate fibration.

Consider now the lifting axiom.  One half of it is immediate from the
definition of a fibration.  For the other half, say $X @>p>> Y$ is an
acyclic intermediate fibration.  Factor $p$ as $X @>j>> Z @>q>>
Y$, where $j$ is in $\Cal C$ and $q$ has the right lifting property
with respect to every map in $\Cal C$.  Then, as before, $q$ is an
acyclic intermediate fibration, so $j$ is acyclic, and by the 
definition of an intermediate fibration, we have a lifting in the square
$$\xymatrix{
X \ar@{=}[r] \ar[d]_j & X \ar[d]^p \\
Z \ar[r]^q \ar@{-->}[ur] & Y
}$$
Hence, $p$ is a retract of $q$, so $p$ also has the right lifting
property with respect to every map in $\Cal C$.

The remaining three axioms for a model structure are clear.  Right 
properness follows from right properness of the injective structure,
and left properness follows from left properness of the projective 
structure.  Finally, Axiom SM7, relating the simplicial structure
and the model structure, may be verified using \cite{GJ, II.3.12}.
\qed \enddemo

Without a useful generating set of acyclic intermediate cofibrations it
is not easy to describe the intermediate fibrations, but for acyclic
intermediate fibrations the following characterization is immediate.

\proclaim{14. Proposition}  An acyclic map $F\to G$ of prestacks is an
intermediate fibration if and only if it has the homotopy right lifting
property with respect to all Stein inclusions.
\endproclaim

\demo{Proof}  By definition of the intermediate structure, an acyclic map
$F\to G$ is an intermediate fibration if and only if there is a lifting in
every square
$$\xymatrix{
S\times\partial\Delta^n\cup_{T\times\partial\Delta^n} T\times\Delta^n
\ar[r] \ar[d]  & F \ar[d] \\
S\times\Delta^n \ar[r] & G
}$$
where $T\hookrightarrow S$ is a Stein inclusion and $n\geq 0$, that is,
by adjunction, in every square
$$\xymatrix{
\partial\Delta^n \ar[r] \ar[d] & F(S) \ar[d] \\
\Delta^n \ar[r] & G(S) \times_{G(T)} F(T)}$$
This means precisely that the map $F(S) \to G(S) 
\times_{G(T)} F(T)$ is an acyclic fibration for
every Stein inclusion $T\hookrightarrow S$.
\qed\enddemo

\subhead 15. The three structures are different  \endsubhead
Two simple examples show that the projective, intermediate, and
injective model structures on $\frak S$ are all different.  First
consider the unit disc $\Bbb D$ (or rather the prestack on $\Cal S$ it
represents).  Since $\Bbb D$ is holomorphically contractible, it is 
projectively fibrant by Theorem 9.  On the other hand, by Liouville's
Theorem, the inclusion $\{0,\frac 1 2\} \hookrightarrow \Bbb D$ does not
factor through the inclusion $\{0,\frac 1 2\} \hookrightarrow \Bbb C$,
which is an intermediate cofibration, so $\Bbb D$ is not intermediately
fibrant.

The complex plane $\Bbb C$ is projectively fibrant for the same reason
that $\Bbb D$ is.  Since $\Bbb C$ is elliptic, it is intermediately fibrant
(see below).  However, $\Bbb C$ is not injectively fibrant; in fact, no 
nondiscrete complex manifold $X$ is.  The inclusion of $\Bbb D$ into the 
disc of radius 2 is a pointwise acyclic monomorphism, but there are many 
holomorphic maps $\Bbb D\to X$ that do not factor through it, so $X$ is not
even coarsely injectively fibrant.

\subhead  16. The Oka property for manifolds and maps  \endsubhead
We say that a holomorphic map $f:X\to Y$ is Oka if it satisfies one of
the following equivalent conditions for every Stein inclusion
$j:T\hookrightarrow S$.

(i) The map $f$ is a topological fibration and satisfies the {\it
Parametric Oka Principle with Interpolation}, meaning that for every
finite polyhedron $P$ with subpolyhedron $Q$ and every diagram
$$\xymatrix{Q \ar[r] \ar[d] & \Cal O(S,X) \ar[r] \ar[d] & \Cal C(S,X)
\ar[d] \\
P \ar[r] \ar@{-->}[urr] \ar@{-->}[ur] & \Cal O(S,Y)\times_{\Cal O(T,Y)}
\Cal O(T,X) \ar[r] & \Cal C(S,Y)\times_{\Cal C(T,Y)} \Cal C(T,X)}$$
of continuous maps,
every lifting $P\to\Cal C(S,X)$ in the big square can be deformed through
liftings in the big square to a lifting that factors through $\Cal O(S,X)$
and is thus a lifting in the left-hand square.  (We recall that a
Serre fibration between smooth manifolds is a Hurewicz fibration \cite{C},
so we will simply call such a map a topological fibration.) 

(ii) A stronger version of condition (i), in which $Q\to P$ is any
cofibration between cofibrant topological spaces and the conclusion is
that the inclusion of the space of liftings $P\to \Cal O(S,X)$ in the
left-hand square into the space of liftings $P\to \Cal C(S,X)$ in the big
square is acyclic.  (Here, and everywhere else in the paper, the notion of
cofibrancy for topological spaces and continuous maps is the stronger one
that goes with Serre fibrations rather than Hurewicz fibrations.)

(iii) The induced map
$$\Cal O(S,X)@>{(f_*,j^*)}>> \Cal O(S,Y)\times_{\Cal O(T,Y)} \Cal O(T,X)$$
is a Serre fibration, and the inclusion
$$\Cal O(S,X)\hookrightarrow \Cal C_{f,T}(S,X) :=
\{h\in\Cal C(S,X) : f\circ h \text{ and } h|T \text{ are holomorphic}\}$$
is acyclic.  Note that $\Cal C_{f,T}(S,X)$ is the pullback of the
right-hand square in condition (i), so when $f$ is a topological
fibration, this inclusion being acyclic is equivalent to that square
being a homotopy pullback.

(iv)  The induced map
$$\Cal O(S,X)@>{(f_*,j^*)}>> \Cal O(S,Y)\times_{\Cal O(T,Y)} \Cal O(T,X)$$
is a Serre fibration, and in any square of holomorphic maps
$$\xymatrix{
T \ar[r] \ar[d] & X \ar[d] \\
S \ar[r] \ar@{-->}[ur] & Y
}$$
the inclusion of the space of holomorphic liftings $S\to X$ into the
space of continuous liftings is acyclic (where these spaces are, as
usual, given the compact-open topology).

\smallskip
Before proving the equivalence of these conditions, we will make a few
remarks.

Observe that the target of $(f_*,j^*)$ is the space of commuting squares
of holomorphic maps in which the map on the left is $j$ and the map on 
the right is $f$.  The fibre over such a square is its set of liftings.
Taking $T=\varnothing$ in each of the conditions gives the weak Oka
property, which we know is equivalent to $f$ being a projective fibration.

Using the Stein inclusion $\varnothing \hookrightarrow \frak p$,
we see that an Oka map is a topological fibration, so its image
is a union of connected components of the target.  An Oka map has
the right lifting property with respect to the inclusion of a point into a
ball, so it is a submersion.  In fact, if a holomorphic map $f:X\to Y$ is
Oka, $q$ is a point in a contractible Stein open subset $V$ of $Y$, and
$p\in f^{-1}(q)$, then condition (iv) implies that $f$ has a holomorphic
section (a right inverse) $V\to X$ taking $q$ to $p$.

A complex manifold $X$ is said to be Oka if the constant map $X\to
\frak p$ is Oka.  This is equivalent to $X$ being weakly Oka and the
restriction map $\Cal O(S,X)\to\Cal O(T,X)$ being a Serre fibration for
every Stein inclusion $T\hookrightarrow S$.  Namely, if $X$ is weakly Oka,
the pullback $\{h\in\Cal C(S,X) : h|T \in\Cal O(T,X)\}\hookrightarrow
\Cal C(S,X)$ of the acyclic map $\Cal O(T,X)\hookrightarrow \Cal C(T,X)$
by the Serre fibration $\Cal C(S,X)\hookrightarrow \Cal C(T,X)$ is
acyclic, and since $\Cal O(S,X)\hookrightarrow \Cal C(S,X)$ is also
acyclic, condition (iii) follows.

\smallskip
Let us now prove the equivalence of the four conditions defining the Oka
property.

(i) $\Rightarrow$ (iii):  Since $f$ is a topological fibration and
$T\hookrightarrow S$ is a topological cofibration, the map $\Cal C(S,X)
\to \Cal C(S,Y)\times_{\Cal C(T,Y)}\Cal C(T,X)$ is a Serre fibration,
so every diagram of continuous maps as below has a lifting as indicated.
$$\xymatrix{[0,1]^n \ar[r] \ar[d] & \Cal O(S,X) \ar[r] \ar[d] & \Cal C(S,X)
\ar[d] \\
[0,1]^{n+1} \ar[r] \ar@{-->}[urr] & \Cal O(S,Y)\times_{\Cal O(T,Y)}
\Cal O(T,X) \ar[r] & \Cal C(S,Y)\times_{\Cal C(T,Y)} \Cal C(T,X)}$$
The Parametric Oka Principle with Interpolation now gives a lifting
$[0,1]^{n+1}\to\Cal O(S,X)$, showing that $\Cal O(S,X)\to \Cal O(S,Y)
\times_{\Cal O(T,Y)} \Cal O(T,X)$ is a Serre fibration.

To prove that $\Cal O(S,X)\hookrightarrow \Cal C_{f,T}(S,X)$ is acyclic,
apply the Parametric Oka Principle with Interpolation to diagrams of the
form
$$\xymatrix{Q \ar[r] \ar[d] & \Cal O(S,X) \ar[r] \ar[d] &
\Cal C_{f,T}(S,X) \ar[dl] \\
P \ar[r] \ar@{-->}[urr] & \Cal O(S,Y)\times_{\Cal O(T,Y)} \Cal O(T,X)}$$
taking $Q\to P$ to be either the inclusion of a point in the $n$-sphere,
$n\geq 1$, or the inclusion of the $n$-sphere in the closed $(n+1)$-ball,
$n\geq -1$.

(iii) $\Leftrightarrow$ (iv):  Consider the diagram
$$\xymatrix{\{\text{holomorphic liftings}\} \ar[r] \ar[d] &
\{\text{continuous liftings}\} \ar[d] \\
\Cal O(S,X) \ar[d] \ar[r] & \Cal C_{f,T}(S,X) \ar[dl] \\
\Cal O(S,Y)\times_{\Cal O(T,Y)} \Cal O(T,X)}$$
The common first part of conditions (iii) and (iv) implies that the
lower downward maps are Serre fibrations, so each horizontal map is
acyclic if and only if the other one is.

(iii) $\Rightarrow$ (ii):  Assume now that $Q\to P$ is any cofibration
between cofibrant topological spaces and consider a diagram as in condition
(i), or equivalently, a diagram
$$\xymatrix{Q \ar[r] \ar[d] & \Cal O(S,X) \ar[r] \ar[d] &
\Cal C_{f,T}(S,X) \ar[d] \\
P \ar[r] & \Cal O(S,Y)\times_{\Cal O(T,Y)}
\Cal O(T,X) \ar@{=}[r] & \Cal O(S,Y)\times_{\Cal O(T,Y)} \Cal O(T,X)}$$
of continuous maps.  Let us write $A=\Cal O(S,X)$, $B=\Cal O(S,Y)
\times_{\Cal O(T,Y)}\Cal O(T,X)$, $C=\Cal C_{f,T}(S,X)$, $\Cal L_A$ for
the space of liftings $P\to A$ in the left-hand square, and $\Cal L_C$
for the space of liftings $P\to C$ in the big square.  Then we have a
diagram
$$\xymatrix{
\Cal L_A \ar[d] \ar[r] & \Cal C(P,A) \ar[r] \ar[d] &
\Cal C(P,B)\times_{\Cal C(Q,B)} \Cal C(Q,A) \ar[d] \\
\Cal L_C \ar[r] & \Cal C(P,C) \ar[r] &
\Cal C(P,B)\times_{\Cal C(Q,B)} \Cal C(Q,C)
}$$
The right-hand horizontal maps are Serre fibrations by Axiom SM7, because
$Q\to P$ is a cofibration and $A\to B$ and $C\to B$ are Serre fibrations.
The middle vertical map is acyclic because $P$ is cofibrant and $A\to C$
is acyclic.  To see that the right-hand vertical map is acyclic,
consider the cube
$$\xymatrix{
\text{pullback} \ar@{-->}[dd] \ar[dr] \ar[rr] & & \Cal C(Q,A) \ar[dd]
\ar[dr] \\
 & \Cal C(P,B) \ar@{=}[dd] \ar[rr] & & \Cal C(Q,B) \ar@{=}[dd] \\
\text{pullback} \ar[dr] \ar[rr] & & \Cal C(Q,C) \ar[dr] \\
 & \Cal C(P,B) \ar[rr] & & \Cal C(Q,B)
}$$
The map $\Cal C(Q,A) \to \Cal C(Q,C)$ is acyclic because $Q$ is
cofibrant and $A\to C$ is acyclic.  The maps $\Cal C(Q,A) \to \Cal C(Q,B)$
and $\Cal C(Q,C) \to \Cal C(Q,B)$ are Serre fibrations because $Q$ is
cofibrant and $A\to B$ and $C\to B$ are Serre fibrations.  Hence, the top
and bottom squares are homotopy pullbacks and we get an induced weak
equivalence of the pullbacks.  Therefore, finally, we get an induced
weak equivalence $\Cal L_A \to \Cal L_C$.

\subhead  17.  Subellipticity and the Oka property  \endsubhead
Subelliptic manifolds satisfy the Parametric Oka Principle with
Interpolation.  This theorem originated in Gromov's work \cite{G} and
was proved in detail by Forstneri\v c and Prezelj \cite{FP, Thm\. 1.4}
for elliptic manifolds; for the extension from ellipticity to 
subellipticity, see \cite{F1}.  Hence, subelliptic manifolds are Oka.

I conjecture that this result extends to nonconstant maps: that a
holomorphic map which is both a subelliptic submersion and a topological
fibration is Oka.  This is an open question even for nontrivial elliptic
fibre bundles.  Here is a small step in this direction, proving the
conjecture in the case of discrete fibres, including the case of covering
maps.

\proclaim{18. Proposition}  A holomorphic map which is a topological
fibration and a local biholomorphism is Oka.
\endproclaim

\demo{Proof}  Let $f:X\to Y$ be a topological fibration and a local
biholomorphism and $T\hookrightarrow S$ be a Stein inclusion.  First,
the inclusion $\Cal O(S,X)\hookrightarrow \Cal C_{f,T}(S,X)$ is acyclic:
it is in fact the identity map because $f$ is a local biholomorphism.
Second, the square
$$\xymatrix{
\Cal O(S,X) \ar[r] \ar[d]_\alpha & \Cal C(S,X) \ar[d]^\beta \\
\Cal O(S,Y) \times_{\Cal O(T,Y)} \Cal O(T,X) \ar[r] & 
\Cal C(S,Y) \times_{\Cal C(T,Y)} \Cal C(T,X)
}$$
is a pullback, because a continuous lifting in a square of holomorphic
maps with right-hand map $f$ is holomorphic, again because $f$ is a local
biholomorphism.  Since $f$ is a Serre fibration, so is $\beta$ by Axiom
SM7, and hence $\alpha$.
\qed\enddemo

We now come to the main result of this paper, describing the intermediate
fibrations.  Notice the similarity with the Oka property as expressed by
condition (iv) above.

\proclaim{19. Theorem (characterization of intermediate fibrations)}
A map $F\to G$ of prestacks is an intermediate fibration if and only if
\roster
\item for every Stein inclusion $T\hookrightarrow S$, the induced map
$$F(S) \to G(S)\times_{G(T)} F(T)$$
is a fibration, and
\item in any diagram
$$\xymatrix{
T \ar[r] \ar[d] & F \ar[d] \ar[r] & PF \ar[d] \\
S \ar[r] & G \ar[r] & PG
}$$
the map of the simplicial set of liftings $S\to F$ into the simplicial
set of liftings $S\to PF$, given by postcomposition with $F\to PF$, is
acyclic.
\endroster
\endproclaim

We remark that taking $T=\varnothing$ in \therosteritem1 and
\therosteritem2 yields precisely the description of projective
fibrations between prestacks fibrant at $\frak p$ given earlier:
\therosteritem1 says that $F\to G$ is a pointwise fibration and
\therosteritem2, using \therosteritem1, says that the induced map
$F\to G\times_{PG}PF$ is pointwise acyclic.

\demo{Proof}  First suppose $F\to G$ is an intermediate fibration and
let $T\hookrightarrow S$ be a Stein inclusion.  Then \therosteritem1
follows directly from Axiom SM7 for the intermediate structure.  As for
\therosteritem2, consider the equivalent diagram
$$\xymatrix{
T \ar[r] \ar[d] & F \ar[d] \ar[r] & E \ar[ld] \\
S \ar[r] & G & 
}$$
where $E=G\times_{PG}PF$.  Since $F(\frak p)\to G(\frak p)$ is a
fibration, $PF\to PG$ and hence $E\to G$ is an injective fibration.
Also, $F\to E$ is acyclic.  Working in the over-under category
$T\downarrow \frak S\downarrow G$ with the model structure induced from
the intermediate structure on $\frak S$, we need to show that
$\Hom_{T\downarrow \frak S\downarrow G}(S,F)\to\Hom_{T\downarrow
\frak S\downarrow G}(S,E)$ is acyclic.  Now $T\to F\to G$ and $T\to E\to
G$ are fibrant in $T\downarrow \frak S\downarrow G$ since $F\to G$ and
$E\to G$ are fibrations in $\frak S$.  By Brown's Lemma \cite{H, 7.7},
we may assume that $F\to E$ is an intermediate fibration in addition to
being acyclic.  Consider the fibration sequences
$$\xymatrix{
\Hom_{T\downarrow \frak S\downarrow G}(S,F) \ar[d] \ar[r] & \Hom(S,F)
\ar[d] \ar[r] & \Hom(S,G)\times_{\Hom(T,G)}\Hom(T,F) \ar[d] \\
\Hom_{T\downarrow \frak S\downarrow G}(S,E) \ar[r] & \Hom(S,E)
\ar[r] & \Hom(S,G)\times_{\Hom(T,G)}\Hom(T,E) }$$
Since $S$ and $T$ are intermediately cofibrant, the middle and
right-hand vertical maps are acyclic, so the left-hand vertical map is
acyclic too.

Now suppose $F\to G$ satisfies \therosteritem1 and \therosteritem2.  We
need to show that $F\to G$ is an intermediate fibration.  First,
\therosteritem1 implies that $\Hom(B,F) \to \Hom(B,G)
\times_{\Hom(A,G)} \Hom(A,F)$ is a fibration for every intermediate
cofibration $A\to B$ (this property is preserved under simplicial
saturation).

Let us say that an intermediate cofibration $A\to B$ is {\it good} if in
any diagram
$$\xymatrix{
A \ar[r] \ar[d] & F \ar[d] \ar[r] & PF \ar[d] \\
B \ar[r] & G \ar[r] & PG
}$$
in $\frak S$, the map of the simplicial set of liftings $B\to F$ into
the simplicial set of liftings $B\to PF$ is acyclic.  By \therosteritem2,
Stein inclusions are good.  With the help of \therosteritem1, we will
show that all intermediate cofibrations are good.  Assuming this,
the proof is complete.  Namely, take a square
$$\xymatrix{
A \ar[r] \ar[d] & F \ar[d] \\
B \ar[r] & G }$$
where $A\to B$ is an acyclic intermediate cofibration.  By
\therosteritem1, $F(\frak p)\to G(\frak p)$ is a fibration, so $PF\to PG$
is an injective fibration and there is a lifting $B\to PF$.  Since
$A\to B$ is good, there is also a lifting $B\to F$.

By piecing together arguments that have already been used in this paper,
the reader can show that the generating cofibrations $S\times\partial
\Delta^n\cup_{T\times\partial\Delta^n}T\times\Delta^n \to S\times\Delta^n$
are good.  It is easy to see that being good is preserved under pushouts
and retracts: pushouts give isomorphisms of lifting spaces and retracts
give retractions of lifting spaces.  It remains to show that a transfinite
composition of good intermediate cofibrations is good.

Let $A\to B$ and $B\to C$ be good and consider a diagram
$$\xymatrix{
A \ar[r] \ar[d] & F \ar[dd] \ar[r] & PF \ar[dd] \\
B \ar[d] & & \\
C \ar[r] & G \ar[r] & PG }$$
Let $\Cal L_{AC}$ and $\Cal L_{AC}'$ be the simplicial sets of liftings
$C\to F$ and $C\to PF$ in the squares with left-hand map $A\to C$ and
right-hand maps $F\to G$ and $PF\to PG$, respectively.  We define
$\Cal L_{AB}$ and $\Cal L_{AB}'$ similarly.  The fibre over a lifting
$B\to F$ of the map $\Cal L_{AC}\to\Cal L_{AB}$ given by precomposing
with $B\to C$ is the simplicial set $\Cal L_{BC}$ of liftings in the
square with left-hand map $B\to C$, right-hand map $F\to G$, and this
particular top map $B\to F$.  We define $\Cal L_{BC}'$ similarly.  We
have a pullback square
$$\xymatrix{
\Cal L_{AC} \ar[r] \ar[d] & \Cal L_{AB} \ar[d] \\
\Hom(C,F) \ar[r] & \Hom(C,G)\times_{\Hom(B,G)}\Hom(B,F) }$$
where the right-hand map takes a map in $\Cal L_{AB}$ to the constant map
$C\to G$ in $\Hom(C,G)$ and itself in $\Hom(B,F)$.  Since the bottom map
is a fibration, so is the top map $\Cal L_{AC} \to \Cal L_{AB}$.  (It
follows that the simplicial set of liftings in any square with
right-hand map $F\to G$ whose left-hand map is an intermediate
cofibration is fibrant: just take $A=B$ and $A\to B$ to be the identity
map.)

By the same argument, $\Cal L_{AC}'\to \Cal L_{AB}'$ is also a
fibration.  Thus the rows in the diagram
$$\xymatrix{
\Cal L_{BC} \ar[r] \ar[d] & \Cal L_{AC} \ar[r] \ar[d] & \Cal L_{AB}
\ar[d] \\
\Cal L_{BC}' \ar[r] & \Cal L_{AC}' \ar[r] & \Cal L_{AB}' }$$
are fibration sequences.  The left-hand and right-hand vertical maps are
acyclic by assumption, so the middle one is too, which shows that
$A\to C$ is good.

We now move to the transfinite case.  Let $\lambda$ be an ordinal and
$A:\lambda\to\frak S$ be a
functor such that for every limit ordinal $\gamma<\lambda$, the induced
map $\colim_{\alpha<\gamma}A_\alpha \to A_\gamma$ is an isomorphism, and
such that for every ordinal $\alpha$ with $\alpha+1<\lambda$, the map
$A_\alpha\to A_{\alpha+1}$ is a good intermediate cofibration.  We will
show by transfinite induction that the composition $A_0\to
\colim_{\alpha<\lambda} A_\alpha$ is good.  Suppose $\mu\leq \lambda$
and $A_0 \to \colim_{\alpha<\beta} A_\alpha$ is good for all
$\beta<\mu$.  We need to show that $A_0\to \colim_{\alpha<\mu}A_\alpha$
is good.

Assume $\mu$ is a successor, say $\mu=\beta+1$.  If $\beta$ is a limit
ordinal, then $A_0\to\colim_{\alpha<\beta} A_\alpha=A_\beta
=\colim_{\alpha<\mu} A_\alpha$ is good by the induction hypothesis.  If
$\beta$ is a successor, say $\beta=\gamma+1$, then $A_0\to\colim_{\alpha
<\beta} A_\alpha=A_\gamma\to A_{\gamma+1}=\colim_{\alpha<\mu} A_\alpha$
is good, being the composition of two good maps. 

Suppose now that $\mu$ is a limit ordinal and take a square
$$\xymatrix{
A_0 \ar[r] \ar[d] & F \ar[d] \\
\colim_{\alpha<\mu}A_\alpha \ar[r] & G }$$
Define a $\mu$-tower $\Cal L:\mu^{\text{op}}\to s\bold{Set}$ such that
$\Cal L_\alpha$ is the simplicial set of liftings in the square
$$\xymatrix{
A_0 \ar[rr]\ar[d] & & F \ar[d] \\
A_\alpha \ar[r] & \colim_{\alpha<\mu}A_\alpha \ar[r] & G }$$
for $\alpha<\mu$.  Define $\Cal L'$ similarly for $PF\to PG$.  Then
$\Cal L$ and $\Cal L'$ are fibrant objects in the category of $\mu$-towers
with the pointwise cofibration simplicial model structure \cite{GJ,
VI.1}, the main point being that for all $\alpha<\mu$, the map $\Cal
L_{\alpha+1}\to \Cal L_\alpha$ is a fibration, as shown above.  Thus,
since the map $\Cal L\to \Cal L'$ is pointwise acyclic by the induction
hypothesis, it induces an acyclic map from $\lim_{\alpha<\mu}
\Cal L_\alpha$ to $\lim_{\alpha<\mu} \Cal L'_\alpha$, that is, from the
simplicial set of liftings $\colim_{\alpha<\mu}A_\alpha\to F$ to the
simplicial set of liftings $\colim_{\alpha<\mu}A_\alpha\to PF$. 
\qed\enddemo

Suppose that the prestacks $F$ and $G$ are represented by complex manifolds
$X$ and $Y$ respectively.  We have
$$\Hom(S,PX)=PX(S)=\Hom(sS,sX)=s\Cal C(|sS|,X).$$
Using the homotopy equivalence $|sS|\to S$, we can verify that our
characterization of the map $F\to G$ induced by a holomorphic map $X\to
Y$ being an intermediate fibration means precisely that $X\to Y$ satisfies
the Oka property as defined by condition (iv) above.

\proclaim{20. Corollary}  A holomorphic map is an intermediate fibration
if and only if it is Oka.
\endproclaim

It follows that subelliptic manifolds are intermediately fibrant and
that holomorphic covering maps are intermediate fibrations.  Also,
the class of Oka maps is preserved under composition, pullbacks, and
retracts.

Our conjecture now looks like this.

\proclaim{21. Conjecture} A subelliptic submersion is an intermediate
fibration if and only if it is a topological fibration.
\endproclaim

This would be a new manifestation of the Oka Principle, in a sense dual
to the usual formulations that refer to Stein manifolds, saying that
for holomorphic maps satisfying the geometric condition of
subellipticity there is only a topological obstruction to being a
fibration in our new, holomorphic sense.

\subhead 22. An alternative approach to the intermediate structure
\endsubhead
We have gone from the coarse projective structure on $\frak S$ to the
intermediate structure via the projective structure by first enlarging
the class of weak equivalences by a Bousfield localization, keeping the
cofibrations fixed, and then enlarging the class of cofibrations,
keeping the weak equivalences fixed.  Alternatively, we could do this
the other way around, passing through what we shall call the coarse
intermediate structure on $\frak S$.  The cofibrations in this structure
are the same as in the intermediate structure, but the weak equivalences
are defined pointwise, and the proof of Theorem 13 goes through word for
word.

A modification of the proof of Theorem 19 gives a characterization
of the coarse intermediate fibrations.  Take a map $F\to G$ of prestacks.
Instead of $PF\to PG$, we now use a coarsely injectively fibrant model
$\tilde F\to \tilde G$ of $F\to G$.  In particular, $F\to\tilde F$
and $G\to\tilde G$ are pointwise acyclic and $\tilde F\to\tilde G$ is a
coarse injective fibration.  Suppose that $F\to G$ satisfies property
\therosteritem1 in Theorem 19.  The key point is that Stein
inclusions are now automatically good, that is, \therosteritem1 implies
\therosteritem2, and we can go on to show that all intermediate
cofibrations are good as before.  Namely, consider a diagram
$$\xymatrix{
T \ar[r] \ar[d] & F \ar[d] \ar[r] & \tilde F \ar[d] \\
S \ar[r] & G \ar[r] & \tilde G
}$$
and the induced diagram of fibration sequences
$$\xymatrix{\{\text{liftings } S\to F\} \ar[r] \ar[d] & F(S) \ar[r]
\ar[d] & G(S)\times_{G(T)}F(T) \ar[d] \\
\{\text{liftings }S\to \tilde F\} \ar[r] & \tilde F(S) \ar[r] &
\tilde G(S)\times_{\tilde G(T)}\tilde F(T)}$$
Since the middle and right-hand vertical arrows are acyclic, so is the
left-hand vertical arrow.  It follows that $F\to G$ is a coarse
intermediate fibration if and only if it satisfies property
\therosteritem1, that is, for every Stein inclusion $T\hookrightarrow S$,
the induced map $F(S) \to G(S)\times_{G(T)}F(T)$ is a fibration.  Hence,
the coarse intermediate fibrations are precisely those maps that have
the right lifting property with respect to the maps
$$S\times\Lambda_k^n\cup_{T\times\Lambda_k^n}
T\times\Delta^n \to S\times\Delta^n,$$
where $\Lambda_k^n$ denotes the $k$-th horn of $\Delta^n$, $0\leq k\leq
n$ (just look at the squares in the proof of Proposition 14).  By a
standard factorization and retraction argument, these maps form a
generating set of acyclic coarse intermediate cofibrations.

Let us now pass to the intermediate structure by a Bousfield
localization.  As we saw for the projective structure earlier, a coarse
intermediate fibration $F\to G$ between prestacks fibrant at $\frak p$
is an intermediate fibration if and only if the square
$$\xymatrix{
F \ar[r] \ar[d] & G \ar[d] \\
PF \ar[r] & PG }$$
is a homotopy pullback in the coarse intermediate structure.  Since $PF\to
PG$ is an injective fibration, this is simply the old condition that the
induced map $F\to G\times_{PG}PF$ be pointwise acyclic, or in other words,
that $F\to G$ be a projective fibration.  This gives one more
characterization of the Oka property, namely condition (iii) above with
$T=\varnothing$ in its second half, which we had previously observed to
be equivalent to (iii) in the case of manifolds.

\subhead 23. Final remarks \endsubhead
We conclude the paper with a few additional words of motivation.  Model
categories are highly nontrivial structures.  Finding them in a new
area of mathematics should be of interest in itself, especially when they
can be shown to be relevant to a topic as deep and important as the Oka
Principle.  The gist of the results in this paper is that analytically
defined Oka properties for complex manifolds and holomorphic maps fit
into a homotopy-theoretic framework in a precise sense: they are
equivalent to fibrancy in suitable model categories containing the
category of complex manifolds.  Our definitions of the Oka property and
the weak Oka property for maps, extending familiar Oka properties of
manifolds, are in fact dictated by abstract homotopy theory.  In short,
we take the point of view that the Oka Principle is about fibrancy. 

It is hoped that this work will eventually have concrete applications
in complex analysis.  Here are three brief remarks in this direction.
First, whether subelliptic submersions that are also topological fibrations
are closed under composition is unknown.  Subelliptic submersions are
not closed under composition and neither is the class of maps with the
property attributed to elliptic bundles by Gromov's theorem (the second
half of condition (iv) above): just consider $\Bbb D\setminus\{0\}
\hookrightarrow \Bbb C \to \frak p$.  Adding to this property a
holomorphic version of Axiom SM7 (the first half of (iv)) yields our Oka
property with all the functorial properties we could wish for.  It
readily implies that if the target is Oka, so is the source, and, if
our conjecture is true, has being a subelliptic submersion and a
topological fibration as a useful geometric sufficient condition.

Second, by the previous section, the Parametric Oka Principle with
Interpolation, as expressed by condition (i), can be verified by
only checking it for acyclic maps $Q \to P$ (giving coarse intermediate
fibrancy) and for $T=\varnothing$ (giving projective fibrancy).
I do not know a direct proof of this (except in the special case of
manifolds, where it is easy).

Third, homotopy theory may shed light on the relationship
between topological and holomorphic contractibility for Stein manifolds.
I believe it is currently unknown whether the former implies the latter.
If we had a suitable weak sufficient condition for coarse intermediate
fibrancy (weaker than subellipticity) satisfied by a topologically
contractible Stein manifold $S$ which did not have the extension
property with respect to some Stein inclusion, then $S$ would not be
intermediately and hence not projectively fibrant and therefore not
holomorphically contractible.  Candidates for such an example exist
in the literature and are being investigated.  The homotopy-theoretic
side of this problem is to distinguish between coarse and fine
intermediate fibrancy for complex manifolds.

\enddocument